\documentclass[12pt]{amsart}
\usepackage{amsmath, amssymb, amsthm, bm, mathtools, enumitem, caption}
\usepackage{hyperref}

\usepackage[noabbrev,capitalize]{cleveref}
\usepackage{adjustbox}
\crefname{equation}{}{}
\usepackage{ytableau}
\usepackage{graphicx, tikz}
\usepackage[textheight=8.95in, textwidth=7in]{geometry}

\usepackage{microtype}

\newtheorem{theorem}{Theorem}[section]
\newtheorem{lemma}[theorem]{Lemma}
\newtheorem{corollary}[theorem]{Corollary}

\newtheorem*{conjecture*}{Conjecture}

\theoremstyle{definition}

\theoremstyle{remark}
\newtheorem*{remark}{Remark}
\newtheorem*{example}{Example}

\numberwithin{equation}{section}

\DeclarePairedDelimiter\floor{\lfloor}{\rfloor}
\DeclarePairedDelimiter\ceil{\lceil}{\rceil}

\newcommand{\N}{\mathbb N}

\newcommand{\fm}{\mathfrak m}

\newcommand{\Q}{\mathbb Q}

\newcommand{\Z}{\mathbb Z}

\newcommand{\lrb}[1]{\left( #1 \right)}

\title[Certain infinite products in terms of MacMahon type series]{Certain infinite products in terms of MacMahon type series}
\date{\today}
\author{Seokho Jin, Badri Vishal Pandey \and Ajit Singh}
\address{Dept. of Mathematics, Chung-Ang University, 84 Heukseok-ro, Dongjak-gu, Seoul 06974, Republic of Korea}
\email{archimed@cau.ac.kr}
\address{University of Cologne, Department of Mathematics and Computer Science,
	Weyertal 86-90,
	50931 Cologne, Germany}
\email{bpandey@uni-koeln.de}
\address{Dept. of Mathematics, University of Virginia, Charlottesville, VA 22904}
\email{ajit18@iitg.ac.in}

\begin{document}
\begin{abstract}Recently, Ono and the third author discovered that the reciprocals of the theta series $(q;q)_\infty^3$ and $(q^2;q^2)_\infty(q;q^2)_\infty^2$ have infinitely many closed formulas in terms of MacMahon's quasimodular forms $A_k(q)$ and $C_k(q)$.
In this article, we use the well-known infinite product identities due to Jacobi, Watson, and  Hirschhorn to derive further such closed formulas for reciprocals of other interesting infinite products. Moreover, with these formulas, we  approximate these reciprocals to arbitrary order simply using MacMahon's functions and {\it MacMahon type} functions. For example, let $\Theta_{6}(q):=\frac{1}{2}\sum_{n\in\Z} \chi_6(n) n q^{\frac{n^2-1}{24}}$ be the theta function corresponding to the odd quadratic character  modulo $6$. Then for any positive integer $n$, we have
$$\frac{1}{\Theta_{6}(q)}= q^{-\frac{3n^2+n}{2}}\sum_{\substack{k=r_1\\ k\equiv n\hspace{-0.2cm}\pmod{2}}}^{r_2}(-1)^{\frac{n-k}{2}}A_{k}(q)C_{\frac{3n-k}{2}}(q)+O(q^{n+1}),$$
where $r_1:=\floor{\frac{3n-1-\sqrt{12n+13}}{3}}+1$ and $r_2:=\ceil{\frac{3n-1+\sqrt{12n+13}}{3}}-1$.
\end{abstract}

\maketitle
\section{Introduction and Statement of Results}

MacMahon, in his seminal work \cite{MacMahon}, studied the following two $q$-series:
\begin{align} 
A_k(q)&:=\sum_{0< t_1<t_2<\cdots<t_k} \frac{q^{t_1+t_2+\cdots+t_k}}{(1-q^{t_1})^2(1-q^{t_2})^2\cdots(1-q^{t_k})^2} \label{Mac1} \\
\intertext{and} 
C_k(q)&:=\sum_{0< t_1<t_2<\cdots<t_k} \frac{q^{2t_1+2t_2+\cdots+2t_k-k}}{(1-q^{2t_1-1})^2(1-q^{2t_2-1})^2\cdots(1-q^{2t_k-1})^2}.\label{Mac2} 
\end{align}
For positive integers $k$, these $q$-series are natural partition generating functions. Namely, we have
\begin{align} 
	A_k(q)&=\sum_{n=1}^\infty \fm(k;n)q^n=\sum_{\substack{0<t_1<t_2<\dots<t_k\\  (m_1,\dots,m_k)\in \N^k}} m_1m_2\cdots m_k q^{m_1t_1+m_2t_2+\dots +m_kt_k} \label{eq:A-k-description} \\
	\intertext{and} 
	C_k(q)&=\sum_{n=1}^\infty \fm_{odd}(k;n)q^n=\sum_{\substack{0<t_1<t_2<\dots<t_k\\  (m_1,\dots,m_k)\in \N^k}} m_1m_2\cdots m_k q^{m_1(2t_1-1)+m_2(2t_2-1)+\dots +m_k(2t_k-1)}. \label{eq:C-k-description}
\end{align}
Furthermore, these series connect these partition functions to disparate areas of mathematics. For example, they arise in the study of  Hilbert schemes,  $q$-multiple zeta-values, representation theory, and  topological string theory (for example, see \cite{R1, R2, R3, R4}). Recently, Andrews and Rose \cite{Andrews-Rose} and Rose \cite{rose} showed that $A_k(q)$ and $C_k(q)$ are quasimodular forms of weight $\leq 2k$ and level 1 and 2 respectively. Since then there has been many interesting research in this direction
(for example, see \cite{AOS, Andrews-Rose, Bachmann, Bringmann, CVO2024, OS, rose}).
\par 
By inspection, one easily sees that for a positive integer $k$, one has
$$A_k(q)=q^{\frac{k(k+1)}{2}}+\cdots\quad \text{and}\quad C_k(q)=q^{k^2}+\cdots,$$
so the order at infinity of these $q$-series as a sequence exhibits quadratic growth. In particular, these series cannot form a natural basis for formal power series in $q$. Nevertheless, this does not rule out the possibility that they could be used to write down closed formulas for specifically exotic or exceptional $q$-series. For any complex number $a$ and $|q|<1$, let $(a;q)_\infty:=\prod_{n=0}^\infty(1-aq^n)$ be the $q$-Pochhammer symbol. Recently in \cite{OS}, Ono and the third author showed that if $k$ is a non-negative integer, then one has
\begin{align}
&\frac{1}{(q;q)_{\infty}^3}=\sum_{n=0}^{\infty}p_3(n)q^{n}=
q^{-\frac{k^2+k}{2}}\sum_{m=k}^\infty \binom{2m+1}{m+k+1}A_{m}(q)\label{OS1}\\
\intertext{and}
&\frac{1}{(q^2;q^2)_{\infty}(q;q^2)_\infty^2}=\sum_{n=0}^{\infty} \overline{p}(n)q^{n}=
q^{-k^2}\sum_{m=k}^\infty \binom{2m}{m+k}C_{m}(q)\label{OS2},
\end{align}
where $p_3(n)$ and $\overline{p}(n)$ denote the number of {\it 3-colored partitions} and the number of {\it overpartitions} of size $n$, respectively.
\par

In view of these results, it is natural to ask whether $A_k(q)$ and $C_k(q)$ enjoy the similar properties with respect to other infinite products extending beyond the reciprocals of $(q;q)_\infty^3$  and $(q^2;q^2)_{\infty}(q;q^2)_\infty^2$. In this paper, we show that this is indeed the case. To this end, we note that
\[
	(q;q)_\infty^3 = 1-3 q+5 q^3-7 q^6+9 q^{10}-11 q^{15}+13 q^{21}-15 q^{28}+17 q^{36}-19 q^{45}+21 q^{55}+\cdots
\]
and
\[
	(q^2;q^2)_{\infty}(q;q^2)_\infty^2= 1-2 q+2 q^4-2 q^9+2 q^{16}-2 q^{25}+2 q^{36}-2 q^{49}+2 q^{64}-2 q^{81}+2 q^{100}+\cdots
\]
are theta functions. As our first result, we consider the infinite product theta series  
\[
	\Theta_{6}(q):=(q;q)_{\infty}^3(q;q^2)_\infty^2=\frac{1}{2}\sum_{n\in\Z} \chi_{6}(n) n q^{\frac{n^2-1}{24}}= 1-5q+7q^2- 11q^5 + 13q^7 - 17q^{12} + 19q^{15} - \cdots,
\]
where 
\[
\chi_{6}(n)=\begin{cases}
						1 & \text{if } n\equiv 1\pmod{6},\\
						-1 & \text{if } n\equiv -1\pmod{6}.
					\end{cases}
\]
Its reciprocal is
$$\sum_{n\ge0} a(n)q^n:=\frac{1}{\Theta_{6}(q)}=1 + 5q + 18q^2 + 55q^3 + 149q^4 + 371q^5 + 867q^6 + \cdots.$$
Using $A_k(q)$ and $C_k(q)$, we give infinitely many closed formulas for $\frac{1}{\Theta_{6}(q)}$. Namely, we have the following theorem.


\begin{theorem}\label{quin}
For any integer $n$, we have
\begin{align*}
	\frac{1}{\Theta_6(q)}
	=(-1)^n q^{-\frac{3n^2+n}{2}}\bigg\{
	&\sum_{\alpha+2\beta=3n+1\atop \alpha,\beta\in\Z}(-1)^\beta\sum_{k_1=|\alpha|}^\infty \binom{2k_1}{k_1+\alpha}A_{k_1}(q)\sum_{k_2=|\beta|}^\infty \binom{2k_2}{k_2+\beta}C_{k_2}(q)\\
	&+\sum_{\alpha+2\beta=3n\atop \alpha,\beta\in\Z}(-1)^\beta\sum_{k_1=|\alpha|}^\infty \binom{2k_1}{k_1+\alpha}A_{k_1}(q)\sum_{k_2=|\beta|}^\infty \binom{2k_2}{k_2+\beta}C_{k_2}(q)\bigg\}.
\end{align*}
\end{theorem}

By Theorem~\ref{quin}, we can compute initial segments of $\frac{1}{\Theta_{6}(q)}$ using a ``small number'' of summands. More precisely, we have the following corollary.
\begin{corollary}\label{Cor-quin} Let $n$ be a positive integer. Then for any non-negative integer $m\le n$, we have
	\[
		a\lrb{m}=\sum_{\substack{k=r_1\\ k\equiv n\hspace{-0.2cm}\pmod{2}}}^{r_2} (-1)^{\frac{3n-k}{2}}\sum_{r=1}^{m-1+\frac{3n^2+n}{2}} \fm(k;r)\fm_{odd}\lrb{\frac{3n-k}{2};m+\frac{3n^2+n}{2}-r},
	\]
		where $r_1:=\floor{\frac{3n-1-\sqrt{12n+13}}{3}}+1$ and $r_2:=\ceil{\frac{3n-1+\sqrt{12n+13}}{3}}-1$.
\end{corollary}	
\begin{example}
	If we take $n=100$ then for $m\le n$, Corollary \ref{Cor-quin}~  gives
	\[
		a\lrb{m}=\sum_{k=45}^{55}(-1)^{k}\sum_{r=1}^{m+15049} \fm(2k;r)\fm_{odd}\lrb{150-k;m+15050-r}.
	\]
We note that only 11 values of $k$ are required to compute $a(m)$.
\end{example}
In view of our last result and \eqref{OS1} and \eqref{OS2}, one might suspect that such results depend on the theory of single variable theta functions. However, this turns out not to be the case. For example, we consider the infinite products  
$$\frac{(q;q)_{\infty}^6}{(q;q^5)_\infty(q^4;q^5)_\infty(q^5;q^5)_\infty}=1-5 q+4 q^2+14 q^3-15 q^4-20 q^5-5 q^6+50 q^7+40 q^8-54 q^9-55 q^{10}+\cdots$$
and
$$\frac{(q;q)_{\infty}^6}{(q^2;q^5)_\infty(q^3;q^5)_\infty(q^5;q^5)_\infty}= 1-6 q+10 q^2+5 q^3-26 q^4+15 q^5-10 q^6+31 q^7+5 q^8-50 q^9+60 q^{10}+\cdots.$$
These $q$-series arise as products of two seemingly independent $q$-series, namely the famous Rogers-Ramanujan functions
\begin{align*}
R_{1,4}(q):=\prod_{n=0}^\infty\frac{1}{(1-q^{5n+1})(1-q^{5n+4})}=\sum_{n=0}^\infty \frac{q^{n^2}}{(1-q)(1-q^2)\cdots (1-q^n)}
\end{align*}
and
\begin{align*}
R_{2,3}(q):=	\prod_{n=0}^\infty\frac{1}{(1-q^{5n+2})(1-q^{5n+3})}=\sum_{n=0}^\infty \frac{q^{n^2+n}}{(1-q)(1-q^2)\cdots (1-q^n)},
\end{align*}
and  the partition generating function studied by Ramanujan \cite[(6.1.4) on page 182]{Andrews-Berndt}
 $$F_5(q):=\sum_{n=0}^\infty q_5(5n) q^n=\frac{(q;q)_{\infty}^6}{(q^5;q^5)_\infty},$$
where $q_m(n)$ denotes the number of $m$-colored partitions of $n$ into an even number of distinct parts minus the number of $m$-colored partitions of $n$ into an odd number of distinct parts.
The reciprocals of $R_{1,4}(q) F_5(q)$ and $R_{2,3}(q) F_5(q)$ have $q$-series
\begin{align*}
	\frac{1}{R_{1,4}(q)F_5(q)}&=1 + 6q + 26q^2 + 91q^3 + 282q^4 + 793q^5 + \cdots =:\sum_{n\ge 0} b_{1,4}(n)q^n \\
	\intertext{and} 
	\frac{1}{R_{2,3}(q) F_5(q)}&=1 + 5q + 21q^2 + 71q^3 + 216q^4 + 597q^5 +\cdots=:\sum_{n\ge 0} b_{2,3}(n)q^n.
\end{align*} 
In the following theorem, we give infinitely many closed formulas for $\frac{1}{R_{1,4}(q)F_5(q)}$ and $\frac{1}{R_{2,3}(q) F_5(q)}$ using convolution sums of $A_k(q)$. Namely,
\begin{theorem}\label{thm:Septuple}
	For any integer $n$, we have
		\begin{align*}
		\frac{1}{R_{1,4}(q)F_5(q)} =q^{-\frac{5n^2-n}{2}}&\bigg\{
		\sum_{\alpha+2\beta=5n+2\atop \alpha,\beta\in\Z}(-1)^\beta\sum_{k_1=|\alpha|}^\infty \binom{2k_1}{k_1+\alpha}A_{k_1}(q)\sum_{k_2=|\beta|}^\infty \binom{2k_2}{k_2+\beta}A_{k_2}(q)\\
		&+\sum_{\alpha+2\beta=5n+1\atop \alpha,\beta\in\Z}(-1)^\beta\sum_{k_1=|\alpha|}^\infty \binom{2k_1}{k_1+\alpha}A_{k_1}(q)\sum_{k_2=|\beta|}^\infty \binom{2k_2}{k_2+\beta}A_{k_2}(q)\\
		&-\sum_{\alpha+2\beta=5n-1\atop \alpha,\beta\in\Z}(-1)^\beta\sum_{k_1=|\alpha|}^\infty \binom{2k_1}{k_1+\alpha}A_{k_1}(q)\sum_{k_2=|\beta|}^\infty \binom{2k_2}{k_2+\beta}A_{k_2}(q)\\
		&-\sum_{\alpha+2\beta=5n-2\atop \alpha,\beta\in\Z}(-1)^\beta\sum_{k_1=|\alpha|}^\infty \binom{2k_1}{k_1+\alpha}A_{k_1}(q)\sum_{k_2=|\beta|}^\infty \binom{2k_2}{k_2+\beta}A_{k_2}(q)\bigg\},
	\end{align*}
	and
	\begin{align*}
		 \frac{1}{R_{2,3}(q) F_5(q)}=q^{-\frac{5n^2-3n}{2}}&\bigg\{
		\sum_{\alpha+2\beta=5n\atop \alpha,\beta\in\Z}(-1)^\beta\sum_{k_1=|\alpha|}^\infty \binom{2k_1}{k_1+\alpha}A_{k_1}(q)\sum_{k_2=|\beta|}^\infty \binom{2k_2}{k_2+\beta}A_{k_2}(q)\\
		&+\sum_{\alpha+2\beta=5n-1\atop \alpha,\beta\in\Z}(-1)^\beta\sum_{k_1=|\alpha|}^\infty \binom{2k_1}{k_1+\alpha}A_{k_1}(q)\sum_{k_2=|\beta|}^\infty \binom{2k_2}{k_2+\beta}A_{k_2}(q)\\
		&-\sum_{\alpha+2\beta=5n-2\atop \alpha,\beta\in\Z}(-1)^\beta\sum_{k_1=|\alpha|}^\infty \binom{2k_1}{k_1+\alpha}A_{k_1}(q)\sum_{k_2=|\beta|}^\infty \binom{2k_2}{k_2+\beta}A_{k_2}(q)\\
		&-\sum_{\alpha+2\beta=5n-3\atop \alpha,\beta\in\Z}(-1)^\beta\sum_{k_1=|\alpha|}^\infty \binom{2k_1}{k_1+\alpha}A_{k_1}(q)\sum_{k_2=|\beta|}^\infty \binom{2k_2}{k_2+\beta}A_{k_2}(q)\bigg\}.
	\end{align*}
\end{theorem} 
By Theorem~\ref{thm:Septuple}, we can compute initial segments of $\frac{1}{R_{1,4}(q)F_5(q)}$ and $\frac{1}{R_{2,3}(q)F_5(q)}$ using a ``small number'' of summands.

\begin{corollary}\label{Cor-Septuple}
	Let $n$ be a positive integer. Then for any non-negative integer $m\le n$, we have
	\begin{multline*}
		b_{1,4}(m) = \sum_{\substack{k=s_1\\k\equiv n-1\pmod{2}}}^{s_2} (-1)^{\frac{n-1-k}{2}} \sum_{v=1}^{m-1+\frac{5n^2-3n}{2}} \fm\lrb{k;v}\fm\lrb{\frac{5n-3-k}{2};m+\frac{5n^2-3n}{2}-v} \\
					 +(2n+1)\sum_{v=1}^{m-1_\frac{5n^2-3n}{2}} \fm\lrb{n;v}\fm\lrb{2n-1;m+\frac{5n^2-3n}{2}-v}
	\end{multline*}
	and
	\begin{multline*}
		b_{2,3}(m)	= \sum_{\substack{k=s_3\\k\equiv n\pmod{2}}}^{s_4} (-1)^{\frac{n-k}{2}} \sum_{v=1}^{m-1+\frac{5n^2-n}{2}} \fm\lrb{k;v}\fm\lrb{\frac{5n-2-k}{2};m+\frac{5n^2-n}{2}-v} \\
					 -(2n-1)\sum_{v=1}^{m-1+\frac{5n^2-n}{2}} \fm\lrb{n-1;v}\fm\lrb{2n;m+\frac{5n^2-n}{2}-v},
	\end{multline*}
	where $s_1:=\floor{\frac{5n-4-\sqrt{40n+41}}{5}}+1$, $s_2:=\ceil{\frac{5n-4+\sqrt{40n+41}}{5}}-1$, $s_3:=\floor{\frac{5n-3-\sqrt{40n+49}}{5}}+1$, and $s_4:=\ceil{\frac{5n-3+\sqrt{40n+49}}{5}}-1$. 
\end{corollary}

\begin{example}
		If we take $n=100$ then for $m\le n$, Corollary \ref{Cor-Septuple}~  gives
	\begin{align*}
		b_{1,4}(m)&= \sum_{k=44}^{56} (-1)^k \sum_{r=1}^{m+24849} \fm\lrb{2k-1;r} \fm\lrb{249-k;m-r+24850} \\
				&\hspace{4cm} +201 \sum_{r=1}^{m+24849} \fm\lrb{100;r} \fm\lrb{199;m-r+24850} \\
		\intertext{and}
		b_{2,3}(m)&=\sum_{k=44}^{56} (-1)^k \sum_{r=1}^{m+24949} \fm\lrb{2k;r} \fm\lrb{249-k;m-r+24950} \\
				&\hspace{4cm} -199 \sum_{r=1}^{m+24949} \fm\lrb{99;r} \fm\lrb{200;m-r+24950}.
	\end{align*}
\end{example}
Theorems \ref{quin} and \ref{thm:Septuple} as well as \eqref{OS1} and \eqref{OS2} give us infinitely many closed formulas in terms of $A_k(q)$ and $C_k(q)$ for particularly exotic or special $q$-series. 
In view of this it is natural to ask whether this phenomenon can hold systematically for infinite families of infinite products. To this end, we show that it is true by revisiting the proofs of \eqref{OS1} and \eqref{OS2} which rely on a beautiful connection between MacMahon functions and certain specializations of Jacobi triple product identity.
For non-negative integers $k$ and $\beta$, we define {\it MacMahon type series} $B_{k,\beta}(q)$ and $D_{k,\beta}(q)$ by
\begin{align}\label{B_k}
	B_{k,\beta}(q):=\sum_{0< m_1<m_2<\cdots<m_k} \frac{q^{m_1+m_2+\cdots+m_k+2k\beta}}{(1-q^{m_1+2\beta})^2(1-q^{m_2+2\beta})^2\cdots(1-q^{m_k+2\beta})^2}
\end{align} 
and 
\begin{align}\label{D_k}
	D_{k,\beta}(q):=\sum_{0< m_1<m_2<\cdots<m_k} \frac{q^{m_1+m_2+\cdots+m_k+k(\beta+1)}}{(1-q^{2m_1+\beta+1})^2(1-q^{2m_2+\beta+1})^2\cdots(1-q^{2m_k+\beta+1})^2}.
\end{align} 
Let $[k]$ denote the set $\{1,2,\dots,k\}$. We define $P([k],n,m)$ by the number of pair of partitions $(\lambda,\mu)$ of $n$ into distinct parts, where the parts of 
$\lambda$ and $\mu$ are both taken from $[k]$ and the difference between the number of parts of  $\lambda$ and $\mu$ is equal to $m$. The  generating function for $P([k],n,m)$ is given by
\begin{align}\label{GF}
	\sum_{m=-k}^k\sum_{n=0}^{k(k+1)}P([k],n,m)q^nz^m:=\prod_{j=1}^k(1+zq^j)(1+z^{-1}q^j).
\end{align}
\begin{theorem}\label{JTP}
	Assuming the above notations, the following are true.

	\noindent
	(1) If $k$ and $\beta$ are non-negative integers, then we have
	\begin{multline*}
		\frac{1}{(q;q)_\infty (q^{2\beta+1};q)_\infty^2} \\
		=q^{-\frac{k(k+1)}{2}-2\beta k-2\beta(2\beta+1)}\bigg\{
		\sum_{m+n=k\atop m\in\{-2\beta,\ldots, 2\beta\},n\in\Z}\sum_{t=0}^{2\beta(2\beta+1)} P([2\beta],t,m)q^t\sum_{s=|n+2\beta|}^\infty \binom{2s}{s+n+2\beta}B_{s,\beta}(q)\\
		+\sum_{m+n=k-1\atop m\in\{-2\beta,\ldots, 2\beta\},n\in\Z}\sum_{t=0}^{2\beta(2\beta+1)} P([2\beta],t,m)q^t\sum_{s=|n+2\beta|}^\infty \binom{2s}{s+n+2\beta}B_{s,\beta}(q)\bigg\}.
	\end{multline*}
	
	\noindent
	(2) If $k$ and $\gamma$ are non-negative integers, then we have
	\begin{multline*}
		\frac{1}{(q^2;q^2)_\infty (q^{\gamma+1};q^2)_\infty^2}\\
		=q^{-k^2-\gamma k-\frac{\gamma(\gamma-1)}{2}}\bigg\{
		\sum_{m+n=k\atop m\in\{-(\gamma-1),\ldots, \gamma-1\},n\in\Z}\sum_{t=0}^{\gamma(\gamma-1)} P([\gamma-1],t,m)q^t\sum_{s=|n+\gamma-1|}^\infty \binom{2s}{s+n+\gamma-1}D_{s,\gamma}(q)\\
		+\sum_{m+n=k-1\atop m\in\{-(\gamma-1),\ldots, (\gamma-1)\},n\in\Z}\sum_{t=0}^{\gamma(\gamma-1)} P([\gamma-1],t,m)q^t\sum_{s=|n+\gamma-1|}^\infty \binom{2s}{s+n+\gamma-1}D_{s,\gamma}(q)\bigg\}.
	\end{multline*}
\end{theorem}
\begin{remark}
Note that as a special case of Theorem \ref{JTP} when $\beta=0$ and $\gamma=0$, we recover Theorem 1.1 and Theorem 1.3 (2) in \cite{OS}, respectively.  
\end{remark}
It is natural to ask whether further MacMahon type series give rise to further results. In this direction, we generalize the MacMahon series by demanding that the summands lie in an arithmetic progression. To this end, let $n\in\N$, then for $\ell\neq n,\frac{n}{2}$ positive integers, we define
\begin{align}\label{eq:A-ell-n-k}
	A_{\ell, n, k}(q):=\sum_{\substack{0<m_1<m_2<\cdots <m_k \\ m_k\equiv \pm \ell \pmod{n}}} \frac{q^{m_1+m_2+\cdots +m_k}}{\lrb{1-q^{m_1}}^2 \lrb{1-q^{m_2}}^2\cdots \lrb{1-q^{m_k}}^2}.
\end{align}
In our next theorem we show that certain family of infinite products have infinitely many formulas in terms of these $q$-series.
\begin{theorem}\label{thm:A-ell-n-k}
	Assuming the above, for all $j\geq 0$, we have
	\begin{align*}
		\prod_{m\geq 1}\frac{\lrb{1+q^{nm}}\lrb{1+q^{2nm-n(j+2)+2\ell}}\lrb{1+q^{2nm+nj-2\ell}}}{\lrb{1-q^{nm}}\lrb{1-q^{nm+\ell-n}}^2\lrb{1-q^{nm-\ell}}^2}
		= q^{-\frac{j}{2}\lrb{jn-2\ell+n}}\sum_{k=|j|}^{\infty} \binom{2k}{j+k} A_{\ell,n,k}(q).
	\end{align*}
\end{theorem}
	\begin{remark}
		We point out that for the cases $\ell=n$ or $\frac{n}{2}$ an analogous result can be proved using similar approach as in the proof of Theorem \ref{thm:A-ell-n-k}. In particular for $n=1,\ell=1$ we recover \eqref{OS1}, and for $n=2,\ell=1$, we recover \eqref{OS2}.
	\end{remark}

\section*{Acknowledgments}
	The authors would like to thank Ken Ono for helpful discussions. The authors also thank Kathrin Bringmann, Kevin Gomez, and Gurinder Singh for their valuable comments on the preliminary draft which improved the writing of the paper. The first author was supported by the National Research Foundation of Korea (NRF) grant funded by the Korea government (MSIT) (No. RS-2023-00253814). The second author has received funding from the European Research Council (ERC) under the European Union’s Horizon 2020 research and innovation programme (grant agreement No. 101001179). The third author thanks the support of a Fulbright Nehru Postdoctoral Fellowship. 
	
\section{Proof of Theorem \ref{quin} and Corollary \ref{Cor-quin}}
\begin{proof}[Proof of Theorem \ref{quin}]
	We recall the quintuple product identity (see (6) of \cite{Gordon}) (see also \cite{Watson})
	\begin{align*}
		\sum_{n=-\infty}^{\infty} q^{\frac{3n^2+n}{2}}(z^{3n}-z^{-3n-1})=(q;q)_\infty\prod_{n=1}^\infty (1-zq^n)(1-z^{-1}q^{n-1})(1-z^{2}q^{2n-1})(1-z^{-2}q^{2n-1}).
	\end{align*}
	
Letting $z\rightarrow -z^2$ and a simple re-indexing gives
	\begin{multline*}
			\sum_{n=-\infty}^{\infty}(-1)^n q^{\frac{3n^2+n}{2}}(z^{6n}+z^{-6n-2})\\
			=(q;q)_\infty(1+z^{-2})\prod_{n=1}^\infty (1+z^2q^n)(1+z^{-2}q^{n})(1-z^{4}q^{2n-1})(1-z^{-4}q^{2n-1}).
	\end{multline*}
	After straightforward algebraic manipulation, we find
	\begin{multline*}
		\sum_{n=-\infty}^{\infty}(-1)^n q^{\frac{3n^2+n}{2}}(z^{6n}+z^{-6n-2})\\
		=(q;q)_\infty(1+z^{-2})\prod_{n=1}^\infty \left((1-q^n)^2+q^{n}(z+z^{-1})^2\right)\left((1-q^{2n-1})^2+q^{2n-1}(i(z^2-z^{-2}))^2\right).
	\end{multline*}
	After factoring out $(q;q)_{\infty}^2(q;q^2)_\infty^2$ from the infinite product,  we obtain
	\begin{multline*}
		\sum_{n=-\infty}^{\infty}(-1)^n q^{\frac{3n^2+n}{2}}(z^{6n}+z^{-6n-2})\\
		=(q;q)_\infty^3(q;q^2)_\infty^2(1+z^{-2})\prod_{n=1}^\infty \left(1+\frac{q^{n}(z+z^{-1})^2}{(1-q^n)^2}\right)\left(1+\frac{q^{2n-1}(i(z^2-z^{-2}))^2}{(1-q^{2n-1})^2}\right).
	\end{multline*}
	Using definitions (\ref{Mac1}) and \eqref{Mac2}, we find that the infinite product on the right, as a convolution of two power series in $(z+z^{-1})^2,$ and $(i(z^2+z^{-2}))^2,$ is the generating function for  MacMahon's series.   
	Therefore, we find that
	\begin{align*}
		&\sum_{n=-\infty}^{\infty}\frac{(-1)^n q^{\frac{3n^2+n}{2}}}{(q;q)_\infty^3(q;q^2)_\infty^2}(z^{6n}+z^{-6n-2})\\
		&\hspace{3cm}=(1+z^{-2})\sum_{k_1=0}^\infty A_{k_1}(q)(z+z^{-1})^{2k_1}\sum_{k_2=0}^\infty C_{k_2}(q)(i(z^2-z^{-2}))^{2k_2}.
	\end{align*}
	Thanks to the Binomial Theorem, followed by a simple shift in the index of summation, and culminating with a change in the order of summation, we obtain
	\begin{align*}
		&\sum_{n=-\infty}^{\infty}\frac{(-1)^n q^{\frac{3n^2+n}{2}}}{(q;q)_\infty^3(q;q^2)_\infty^2}(z^{6n}+z^{-6n-2})\\
		&\hspace{1cm}=\left(1+z^{-2}\right)\sum_{k_1=0}^{\infty}A_{k_1}(q)\sum_{r=0}^{2k_1}\binom{2k_1}{r}z^{2r-2k_1}\sum_{k_2=0}^{\infty}C_{k_2}(q)\sum_{s=0}^{2k_2}(-1)^{k_2+s}\binom{2k_2}{s}z^{4s-4k_2}\\
		&\hspace{1cm}=\left(1+z^{-2}\right)\sum_{r=-\infty}^{\infty}\sum_{k_1=|r|}^{\infty}\binom{2k_1}{r+k_1}A_{k_1}(q)z^{2r}\sum_{s=-\infty}^{\infty}\sum_{k_2=|s|}^{\infty}(-1)^s\binom{2k_2}{s+k_2}C_{k_2}(q)z^{4s}.
	\end{align*}
	After multiplying through $(1+z^{-2})$, we obtain
	\begin{multline*}
		\sum_{n=-\infty}^{\infty}\frac{(-1)^n q^{\frac{3n^2+n}{2}}}{(q;q)_\infty^3(q;q^2)_\infty^2}(z^{6n}+z^{-6n-2})\\
		=\sum_{k=-\infty}^{\infty}\left(\sum_{r+2s=k \atop r,s\in\Z}(-1)^s\sum_{k_1=|r|}^{\infty}\binom{2k_1}{r+k_1}A_{k_1}(q)\sum_{k_2=|s|}^{\infty}\binom{2k_2}{s+k_2}C_{k_2}(q)\right.\\
		\left.+ \sum_{r+2s=k+1 \atop r,s\in\Z}(-1)^s\sum_{k_1=|r|}^{\infty}\binom{2k_1}{r+k_1}A_{k_1}(q)\sum_{k_2=|s|}^{\infty}\binom{2k_2}{s+k_2}C_{k_2}(q)\right)z^{2k}.
	\end{multline*}
	The theorem follows by comparing the coefficient of $z^{6k}$ on both sides.
\end{proof}
 
In order to prove Corollary \ref{Cor-quin}, we first prove the following lemma.
\begin{lemma}\label{lem:quin} For any positive integer $n$, we have
	$$\frac{1}{\Theta_6(q)}= q^{-\frac{3n^2+n}{2}}\sum_{\substack{k=r_1\\ k\equiv n\hspace{-0.2cm}\pmod{2}}}^{r_2}(-1)^{\frac{3n-k}{2}}A_{k}(q)C_{\frac{3n-k}{2}}(q)+O(q^{n+1}),$$
	where $r_1$ and $r_2$ are defined in Corollary \ref{Cor-quin}.
\end{lemma}
\begin{proof}
	To prove the result, we want to gather all parts from the right hand side of Theorem \ref{quin} which contribute to the main term. From both the outer summations in the right hand side, we want the pairs $\lrb{\alpha,\beta}\in\Z^2$ with $\alpha+2\beta=3n$, or $3n+1$ which satisfy
	\begin{align}\label{eq:K-apha-beta-ineq}
		K(\alpha,\beta):=\frac{\alpha^2+|\alpha|}{2} +\beta^2 < \frac{3n^2+n}{2}+n+1=3 \cdot\frac{n^2+n}{2}+1,
	\end{align}
	and from the inner sums, we want pairs $(k_1,k_2)$ such that $k_1\ge |\alpha|$ and $k_2\ge |\beta|$ with the condition that
	\begin{align*}
		K(k_1,k_2)< 3 \cdot\frac{n^2+n}{2}+1.
	\end{align*}
	Case 1.1) $\alpha+2\beta=3n+1$ with $\alpha\ge 0$: we have
	\begin{align*}
		K\lrb{\alpha,\frac{3n+1-\alpha}{2}} = \frac{3}{4}\alpha^2 -\frac{3n\alpha}{2} + \frac{1}{4}\lrb{9n^2+6n+1}.
	\end{align*}
	This function is minimized when $\alpha=n$, but in this case $\beta=n+\frac{1}{2}\not\in \Z,$ and $K\lrb{n,n+\frac{1}{2}} = 3\cdot \frac{n^2+n}{2}+\frac{1}{4}$, and at all the other values \eqref{eq:K-apha-beta-ineq} is not satisfied.
	
	\noindent
	Case 1.2) $\alpha+2\beta=3n+1$ with $\alpha<0$: a similar argument as in the previous case. gives us that there are no values that satisfy \eqref{eq:K-apha-beta-ineq}.
	
\noindent
	Case 2.1) $\alpha+2\beta=3n$ with $\alpha\ge 0$: we have
	\begin{align*}
		K\lrb{\alpha,\frac{3n-\alpha}{2}} = \frac{3}{4}\alpha^2 +\frac{\lrb{1-3n}\alpha}{2} + \frac{9n^2}{4}.
	\end{align*}
	This function is minimized for integer values $\alpha=n=\beta$, and the minimum value is $\frac{3n^2+n}{2}$, which satisfies \eqref{eq:K-apha-beta-ineq}. We now check for other feasible values. For that we solve \eqref{eq:K-apha-beta-ineq}, i.e., we have
	\begin{align*}
		\frac{3}{4}\alpha^2 +\frac{\lrb{1-3n}\alpha}{2} + \frac{9n^2}{4} &< \frac{3n^2+n}{2}+n+1=3 \cdot\frac{n^2+n}{2}+1,
	\end{align*}
	which implies the condition
	\begin{align}\label{eq:alpha-solution}
	-1\le \frac{3n-1-\sqrt{12n+13}}{3} < \alpha < \frac{3n-1+\sqrt{12n+13}}{3}.
	\end{align}

	\noindent
	Case 2.2) $\alpha+2\beta=3n$ with $\alpha<0$: we have 
	\begin{align*}
		K\lrb{\alpha,\beta} = \frac{\alpha^2-\alpha}{2} +\beta^2.
	\end{align*}
	One easily checks that no solution exists in this case.
	
	Finally, for the inner summations, we have $k_1\ge |\alpha|$ and $k_2\ge |\beta|$, which imply $k_1+2k_2\ge 3n$. We want to make sure that 
	\begin{equation}\label{eq:k_1-k_2}
	\frac{k_1^2+k_1}{2} +k_2^2< 3\cdot\frac{n^2+n}{2}+1.
	\end{equation}
	Since the case of $k_1+2k_2= 3n$ is already covered (when $k_1=|\alpha|$ and $k_2=\beta$), we can assume $k_1+2k_2\ge 3n+1$, whereby we have
	\begin{align*}
		\frac{k_1^2+k_1}{2} +k_2^2 \ge \frac{k_1^2+k_1}{2} +\lrb{\frac{3n+1-k_1}{2}}^2 \ge 3\cdot\frac{n^2+n}{2} + \frac{1}{4}.
	\end{align*}
	Using this one can check no integer pairs $\lrb{k_1,k_2}$ satisfy \eqref{eq:k_1-k_2}. Hence, all the solutions come from \eqref{eq:alpha-solution}. This completes the proof of the lemma.
\end{proof}
\begin{proof}[Proof of Corollary \ref{Cor-quin}]
	Putting the $q$-series expansions \eqref{eq:A-k-description} and \eqref{eq:C-k-description} in Lemma \ref{lem:quin} and comparing the appropriate coefficients completes the proof.
\end{proof}

\section{Proof of Theorem \ref{thm:Septuple} and Corollary \ref{Cor-Septuple}}
\begin{proof}[Proof of Theorem \ref{thm:Septuple}]
	We first recall the septuple product identity from \cite{Hirschhorn1}
	\begin{multline*}
		(q^2;q^5)_\infty(q^3;q^5)_\infty(q^5;q^5)_\infty\sum_{n=-\infty}^{\infty} (-1)^nq^{\frac{5n^2}{2}}\left(q^{\frac{-3n}{2}}z^{5n}+q^{\frac{3n}{2}}z^{5n+3}\right)\\
		-(q;q^5)_\infty(q^4;q^5)_\infty(q^5;q^5)_\infty\sum_{n=-\infty}^{\infty} (-1)^nq^{\frac{5n^2}{2}}\left(q^{\frac{-n}{2}}z^{5n+1}+q^{\frac{n}{2}}z^{5n+2}\right)\\
	=(q;q)_\infty^2\prod_{n=0}^\infty (1-zq^{n})(1-z^{-1}q^{n+1})(1-z^{2}q^{n})(1-z^{-2}q^{n+1}).
	\end{multline*}
Letting $z\rightarrow -z^2$ and doing straightforward algebraic manipulation, we find
	\begin{multline*}
		(q^2;q^5)_\infty(q^3;q^5)_\infty(q^5;q^5)_\infty\sum_{n=-\infty}^{\infty} q^{\frac{5n^2}{2}}\left(q^{\frac{-3n}{2}}z^{10n}-q^{\frac{3n}{2}}z^{10n+6}\right)\\
		-(q;q^5)_\infty(q^4;q^5)_\infty(q^5;q^5)_\infty\sum_{n=-\infty}^{\infty} q^{\frac{5n^2}{2}}\left(q^{\frac{n}{2}}z^{10n+4}-q^{\frac{-n}{2}}z^{10n+2}\right)\\
		=(q;q)_\infty^2(1+z^2)(1-z^4)\prod_{n=1}^\infty \left((1-q^n)^2+q^{n}(z+z^{-1})^2\right)\left((1-q^{n})^2+q^{n}(i(z^2-z^{-2}))^2\right).
	\end{multline*}
	After factoring out $(q;q)_{\infty}^4$ from the infinite product,  we obtain
	\begin{multline*}
		(q^2;q^5)_\infty(q^3;q^5)_\infty(q^5;q^5)_\infty\sum_{n=-\infty}^{\infty} q^{\frac{5n^2}{2}}\left(q^{\frac{-3n}{2}}z^{10n}-q^{\frac{3n}{2}}z^{10n+6}\right)\\
		-(q;q^5)_\infty(q^4;q^5)_\infty(q^5;q^5)_\infty\sum_{n=-\infty}^{\infty} q^{\frac{5n^2}{2}}\left(q^{\frac{n}{2}}z^{10n+4}-q^{\frac{-n}{2}}z^{10n+2}\right)\\
		=(q;q)_\infty^6(1+z^2)(1-z^4)\prod_{n=1}^\infty \left(1+\frac{q^{n}(z+z^{-1})^2}{(1-q^n)^2}\right)\left(1+\frac{q^{n}(i(z^2-z^{-2}))^2}{(1-q^{n})^2}\right).
	\end{multline*}
	Using definition (\ref{Mac1}), we find that the infinite product on the right, as a convolution of two power series in $(z+z^{-1})^2,$ and $(i(z^2+z^{-2}))^2,$ is the generating function for  MacMahon's series.   
	Therefore, we find that
	\begin{multline*}
		\frac{(q^2;q^5)_\infty(q^3;q^5)_\infty(q^5;q^5)_\infty}{(q;q)_\infty^6}\sum_{n=-\infty}^{\infty} q^{\frac{5n^2}{2}}\left(q^{\frac{-3n}{2}}z^{10n}-q^{\frac{3n}{2}}z^{10n+6}\right)\\
		-\frac{(q;q^5)_\infty(q^4;q^5)_\infty(q^5;q^5)_\infty}{(q;q)_\infty^6}\sum_{n=-\infty}^{\infty} q^{\frac{5n^2}{2}}\left(q^{\frac{n}{2}}z^{10n+4}-q^{\frac{-n}{2}}z^{10n+2}\right)\\
		=(1+z^2)(1-z^4)\sum_{k_1=0}^\infty A_{k_1}(q)(z+z^{-1})^{2k_1}\sum_{k_2=0}^\infty A_{k_2}(q)(i(z^2-z^{-2}))^{2k_2}.
	\end{multline*}
	Employing the Binomial Theorem, followed by a simple shift in the index of summation, and finally changing the order of summation, we obtain
	\begin{multline*}
		\frac{(q^2;q^5)_\infty(q^3;q^5)_\infty(q^5;q^5)_\infty}{(q;q)_\infty^6}\sum_{n=-\infty}^{\infty} q^{\frac{5n^2}{2}}\left(q^{\frac{-3n}{2}}z^{10n}-q^{\frac{3n}{2}}z^{10n+6}\right)\\
		-\frac{(q;q^5)_\infty(q^4;q^5)_\infty(q^5;q^5)_\infty}{(q;q)_\infty^6}\sum_{n=-\infty}^{\infty} q^{\frac{5n^2}{2}}\left(q^{\frac{n}{2}}z^{10n+4}-q^{\frac{-n}{2}}z^{10n+2}\right)\\
		=(1+z^2)(1-z^4)\sum_{r=-\infty}^{\infty}\sum_{k_1=|r|}^{\infty}\binom{2k_1}{r+k_1}A_{k_1}(q)z^{2r}\sum_{s=-\infty}^{\infty}\sum_{k_2=|s|}^{\infty}(-1)^s\binom{2k_2}{s+k_2}A_{k_2}(q)z^{4s}.
	\end{multline*}
	After multiplying through $(1+z^2)(1-z^4)$, we obtain
	\begin{align*}
		&\frac{(q^2;q^5)_\infty(q^3;q^5)_\infty(q^5;q^5)_\infty}{(q;q)_\infty^6}\sum_{n=-\infty}^{\infty} q^{\frac{5n^2}{2}}\left(q^{\frac{-3n}{2}}z^{10n}-q^{\frac{3n}{2}}z^{10n+6}\right)\\
		&\hspace{1cm}-\frac{(q;q^5)_\infty(q^4;q^5)_\infty(q^5;q^5)_\infty}{(q;q)_\infty^6}\sum_{n=-\infty}^{\infty} q^{\frac{5n^2}{2}}\left(q^{\frac{n}{2}}z^{10n+4}-q^{\frac{-n}{2}}z^{10n+2}\right)\\
		&\hspace{2cm}=\sum_{k=-\infty}^{\infty}\left(\sum_{r+2s=k \atop r,s\in\Z}(-1)^s\sum_{k_1=|r|}^{\infty}\binom{2k_1}{r+k_1}A_{k_1}(q)\sum_{k_2=|s|}^{\infty}\binom{2k_2}{s+k_2}A_{k_2}(q)\right.\\
		&\hspace{3cm}+\left.\sum_{r+2s=k-1 \atop r,s\in\Z}(-1)^s\sum_{k_1=|r|}^{\infty}\binom{2k_1}{r+k_1}A_{k_1}(q)\sum_{k_2=|s|}^{\infty}\binom{2k_2}{s+k_2}A_{k_2}(q)\right.\\
		&\hspace{3cm}-\left.\sum_{r+2s=k-2 \atop r,s\in\Z}(-1)^s\sum_{k_1=|r|}^{\infty}\binom{2k_1}{r+k_1}A_{k_1}(q)\sum_{k_2=|s|}^{\infty}\binom{2k_2}{s+k_2}A_{k_2}(q)\right.\\
		&\hspace{3cm}\left.- \sum_{r+2s=k-3 \atop r,s\in\Z}(-1)^s\sum_{k_1=|r|}^{\infty}\binom{2k_1}{r+k_1}A_{k_1}(q)\sum_{k_2=|s|}^{\infty}\binom{2k_2}{s+k_2}A_{k_2}(q)\right)z^{2k}.
	\end{align*}
	We complete the proof of the theorem by comparing the coefficient of $z^{10n+4}$ and $z^{10n}$ on both sides.
\end{proof}

We now prove Corollary \ref{Cor-Septuple}. For that, we first prove the following lemma.
\begin{lemma}\label{lem:sept}
For any integer $n\geq 2$, we have
	\begin{multline*}
	\frac{1}{R_{1,4}(q)F_5(q)}\\
= q^{-\frac{5n^2-3n}{2}}\left\{\sum_{\substack{k=s_1\\ k\equiv n-1\hspace{-0.2cm}\pmod{2}}}^{s_2}(-1)^{\frac{n-1-k}{2}}A_{k}(q)A_{\frac{5n-3-k}{2}}(q)+(2n+1)A_n(q)A_{2n-1}(q) \right\}+O(q^{n+1})
\end{multline*}
and 
	\begin{multline*}
	\frac{1}{R_{2,3}(q)F_5(q)}\\
=q^{-\frac{5n^2-n}{2}}\left\{\sum_{\substack{k=s_3\\ k\equiv n\hspace{-0.2cm}\pmod{2}}}^{s_4}(-1)^{\frac{n-k}{2}}A_{k}(q)A_{\frac{5n-2-k}{2}}(q)-(2n-1)A_{n-1}(q)A_{2n}(q)\right\}+O(q^{n+1}),
	\end{multline*}
where $s_1, s_2, s_3,$ and $s_4$ are defined in Corollary \ref{Cor-Septuple}.
\end{lemma}

\begin{proof}
We give a proof for the first part only, since the other part can be obtained in a similar fashion. Recalling that $A_k(q)=q^{\frac{k(k+1)}{2}}(1+O(q))$, as we did in the proof of Corollary \ref{Cor-quin}, what we need is to pick up 
the pairs $(\alpha,\beta)$ such that $K(\alpha,\beta)=\frac{|\alpha|(|\alpha|+1)}{2}+\frac{|\beta|(|\beta|+1)}{2}$ is less than $\frac{5n^2-3n}{2}+n+1=\frac{5n^2-n}{2}+1$.  

\noindent Case 1) $\alpha\geq 0, \beta\geq 0$: In this case $K(\alpha,\beta)=\frac{\alpha^2+\alpha+\beta^2+\beta}{2}=\frac{5\alpha^2-2\alpha(\gamma-1)+\gamma(\gamma+2)}{8}$. Note that for $\gamma=\alpha+2\beta= 5n, 5n-1, 5n-2, 5n-3$, the minimum of $K(\alpha,\beta)$ is attained at $(\alpha, \gamma)=(n-1, 5n-3)$ with the value $\frac{5n^2-3n}{2}$ for integers $\alpha$. We thus see that the only possible cases where the $q$-exponent of the initial term is less than $\frac{5n^2-n}{2}+1$ are when $\gamma=5n-2$ and $\gamma=5n-3$, where the minimizing $\alpha$ is $n-\frac{3}{5}$ and $n-\frac{4}{5}$, respectively. We solve the inequalty $K(\alpha,\beta)<\frac{5n^2-n}{2}+1$ in each case. When $\gamma=5n-2$, we have $$K(\alpha,\beta)=\frac{5\alpha^2-2\alpha(5n-3)+(5n-2)(5n)}{8}=\frac{5}{8}(\alpha+\frac{3-5n}{5})^2+\frac{5n^2-n}{2}-\frac{9}{40}<\frac{5n^2-n}{2}+1,$$ from which we have $\alpha= n-1, n$. But $\alpha=n-1$ is impossible since then $\beta$ is not integral, so in this case we have one solution $(\alpha,\beta)=(n, 2n-1)$. For the case $\gamma=5n-3$, we have $$K(\alpha,\beta)=\frac{5\alpha^2-2\alpha(5n-4)+(5n-3)(5n-1)}{8}=\frac{5}{8}(\alpha+\frac{4-5n}{5})^2+\frac{5n^2-3n}{2}-\frac{1}{40}<\frac{5n^2-n}{2}+1,$$ from which we have $\frac{5n-4-\sqrt{40n+41}}{5}<\alpha<\frac{5n-4+\sqrt{40n+41}}{5}$. Note also that $\frac{5n-4-\sqrt{40n+41}}{5}>-1$ for $n\geq 2$, so that $\floor{\frac{5n-4-\sqrt{40n+41}}{5}}+1\geq 0$. 

For terms other than the initial terms, we need to find $(k_1, k_2)$ such that $k_1\geq |\alpha|$, $k_2\geq |\beta|$ and $\frac{k_1^2+k_1}{2}+\frac{k_2^2+k_2}{2}<\frac{5n^2-n}{2}+1$. First when $\gamma=\alpha+2\beta=5n-2$, the cases $(k_1, k_2)=(n+1, 2n-1), (n, 2n)$ do not satisfy the condition, hence there are no other terms than the initial term $(n, 2n-1)$. When $\gamma=5n-3$, the other term $(k_1, k_2)$ that are possible is when $(k_1, k_2)=(n, 2n-1)$ where $k_1+2k_2=5n-2$. 

Therefore, in this case we have found $(\alpha, \beta)=(n,2n-1)$ or $(\alpha, \frac{5n-3-\alpha}{2})$ such that $\frac{5n-4-\sqrt{40n+41}}{5}<\alpha<\frac{5n-4+\sqrt{40n+41}}{5}$, and $(k_1, k_2)=(n, 2n-1)$ where $\alpha+2\beta=5n-3$.

\noindent Case 2) $\alpha\geq 0, \beta< 0$: In this case $K(\alpha,\beta)=\frac{\alpha^2+\alpha+\beta^2-\beta}{2}=\frac{5\alpha^2-2\alpha(\gamma-3)+\gamma(\gamma-2)}{8}$. By computing the minimums of $K(\alpha,\beta)$ for $\gamma= 5n, 5n-1, 5n-2$, and $5n-3$, we see that in every case the $q$-exponent of the initial term is less than $\frac{5n^2-n}{2}+1$, with minimizing $\alpha=\frac{\gamma-3}{5}$ for each case. However, in this case the condition $\beta<0$ forces $\alpha$ to be larger than $\gamma$, so we must consider the inequality $K(\alpha,\beta)<\frac{5n^2-n}{2}+1$ in the region $[\gamma+1, \infty)$, where the function $K(\alpha,\beta)$ is strictly increasing. But computation shows that $K(\gamma+1,\beta)$ is minimal when $\gamma=5n-3$, with value $\frac{25}{2}n^2-\frac{15}{2}n+\frac{11}{8}$ greater than $\frac{5n^2-n}{2}+1$ if $n\geq 1$. Therefore there is no $\alpha$ in this case.

\noindent Case 3) $\alpha<0$ (hence $\beta=\frac{\gamma-\alpha}{2}>0$): In this case $K(\alpha,\beta)=\frac{\alpha^2-\alpha+\beta^2-\beta}{2}=\frac{5\alpha^2-2\alpha(\gamma+3)+\gamma(\gamma+2)}{8}$. By computing the minimums of $K(\alpha,\beta)$ for $\gamma= 5n, 5n-1, 5n-2$, and $5n-3$, we see that the only possible cases where the $q$-exponent of the initial term is less than $\frac{5n^2-n}{2}+1$ are when $\gamma=5n-2$ and $\gamma=5n-3$, where the minimizing $\alpha$ is $n+\frac{1}{5}$ and $n$ respectively. But since $\alpha<0$ we must consider the inequality $K(\alpha,\beta)<\frac{5n^2-n}{2}+1$ in the region $(-\infty,-1]$, where the function $K(\alpha,\beta)$ is strictly decreasing. But computation shows that $K(-1,\beta)$ is minimal when $\gamma=5n-3$, with value $\frac{25}{2}n^2-\frac{5}{4}n+1$ greater than $\frac{5n^2-n}{2}+1$ if $n\geq 2$. This shows that there is no $\alpha$ in this case. 
\end{proof}
\begin{proof}[Proof of Corollary \ref{Cor-Septuple}]
	Putting the $q$-series expansion from \eqref{eq:A-k-description} in Lemma \ref{lem:sept} and comparing the appropriate coefficients completes the proof.
\end{proof}

\section{Proof of Theorem \ref{JTP} and Theorem \ref{thm:A-ell-n-k}}
Here, we recall a few definitions and earlier work of Rose \cite{rose} that will be useful in the proof of  Theorem \ref{thm:A-ell-n-k}. We first define the Jacobi-theta function.
For $r\in\Q$, we define
		\begin{align*}
			\vartheta_{r}(z;q) := \sum_{m\in\Z+r} z^m q^{\frac{1}{2}m^2}.
		\end{align*}
		As a special case, letting $z\rightarrow -1$, we have
		\begin{align*}
			\vartheta_r(q):=\sum_{m\in\Z+r} (-1)^m q^{\frac{1}{2}m^2}.
		\end{align*}
We also need a special form of the Jacobi triple product (for $r\neq \frac{1}{2}$) \cite[(5)]{rose}
\begin{align}\label{eq:JPR-special-1}
	\vartheta_{r}(z;q) = z^r q^{\frac{r^2}{2}} \prod_{m>0} \lrb{1-q^m}\lrb{1+zq^{m+r-\frac{1}{2}}}\lrb{1+z^{-1}q^{m-r-\frac{1}{2}}}.
\end{align}
Rose proved that the function $A_{\ell,n,k}(q)$ is a quasi-modular form as a special case of \cite[Theorem~1.9]{rose}. We write down this special case as we require it in the proof of Theorem \ref{thm:A-ell-n-k}.
\begin{theorem}\label{thm:rose}
	For $n\in\N,$ and $\ell\neq n,\frac{n}{2}$, we have
	\begin{align*}
		\sum_{k\ge 0} (-1)^k A_{\ell,n,k}(q) x^{2k} = \dfrac{\vartheta_{\frac{\ell}{n}-\frac{1}{2}}(q^n;-z)\vartheta_{-\frac{\ell}{n}+\frac{1}{2}}(q^n;-z)}{\vartheta_{\frac{\ell}{n}-\frac{1}{2}}(q^n)\vartheta_{-\frac{\ell}{n}+\frac{1}{2}}(q^n)},
	\end{align*}
	where $x=z^{\frac{1}{2}}-z^{-\frac{1}{2}}$.
\end{theorem}
\begin{proof}[Proof of Theorem~\ref{thm:A-ell-n-k}]
	From Theorem \ref{thm:rose}, for symmetric set $S=\{\ell,n-\ell\}$ with respect to $n$, we have
	\begin{align} \label{eq:sum-A-ell-n-k}
		\sum_{k=0}^\infty (-1)^k A_{\ell, n,k}(q) \lrb{z^{\frac{1}{2}}-z^{-\frac{1}{2}}}^{2k} &= \dfrac{\vartheta_{\frac{\ell}{n}-\frac{1}{2}}(q^n;-z)\vartheta_{-\frac{\ell}{n}+\frac{1}{2}}(q^n;-z)}{\vartheta_{\frac{\ell}{n}-\frac{1}{2}}(q^n)\vartheta_{-\frac{\ell}{n}+\frac{1}{2}}(q^n)}.
	\end{align}
	Using \eqref{eq:JPR-special-1} the denominator of the right-hand side becomes
	\begin{align}\label{eq:donom-of-RHS}
		\vartheta_{\frac{\ell}{n}-\frac{1}{2}}(q^n)\vartheta_{-\frac{\ell}{n}+\frac{1}{2}}(q^n) = q^{n\lrb{\frac{\ell}{n}-\frac{1}{2}}^2}\prod_{m>0}\lrb{1-q^{nm}}^2\lrb{1-q^{nm+\ell-n}}^2\lrb{1-q^{nm-\ell}}^2.
	\end{align}
	Now, we look at the numerator of the right-hand side of \eqref{eq:sum-A-ell-n-k}. We have
	\begin{align}
		\vartheta_{\frac{\ell}{n}-\frac{1}{2}}(q^n;-z)\vartheta_{-\frac{\ell}{n}+\frac{1}{2}}(q^n;-z) &= \sum_{m\in\Z+\frac{\ell}{n}-\frac{1}{2}} q^{\frac{nm^2}{2}} (-z)^m \sum_{r\in\Z-\frac{\ell}{n}+\frac{1}{2}}q^{\frac{nr^2}{2}} (-z)^r \nonumber\\
					&= \sum_{t\in\Z} (-z)^{t} q^{\frac{n}{2}\lrb{t^2-2t\lrb{\frac{\ell}{n}-\frac{1}{2}}+2\lrb{\frac{\ell}{n}-\frac{1}{2}}^2}} \sum_{m\in\Z}(q^{2n})^{\frac{m^2}{2}} \lrb{q^{-n\lrb{t-\frac{2\ell}{n}+1}}}^m.\nonumber 
	\end{align}
					Using \eqref{eq:JPR-special-1} in the inner sum, we get 
		\begin{align}
					\vartheta_{\frac{\ell}{n}-\frac{1}{2}}(q^n;-z)\vartheta_{-\frac{\ell}{n}+\frac{1}{2}}(q^n;-z) &=\sum_{t\in\Z} (-z)^{t} q^{\frac{n}{2}\lrb{t^2-2t\lrb{\frac{\ell}{n}-\frac{1}{2}}}}q^{n\lrb{\frac{\ell}{n}-\frac{1}{2}}^2} \nonumber\\
									&\hspace{1cm} \prod_{m>0} \lrb{1-q^{2nm}}\lrb{1+q^{2nm-n(t+2)+2\ell}}\lrb{1+q^{2nm+nt-2\ell}}.	\label{eq:num-of-RHS}	
	\end{align}
	We put \eqref{eq:donom-of-RHS} and \eqref{eq:num-of-RHS} back into \eqref{eq:sum-A-ell-n-k} and obtain
	\begin{align*}
		\sum_{t\in\Z} (-z)^{t} q^{\frac{n}{2}\lrb{t^2-2t\lrb{\frac{\ell}{n}-\frac{1}{2}}}} &\prod_{m>0} \frac{\lrb{1-q^{2nm}}\lrb{1+q^{2nm-n(t+2)+2\ell}}\lrb{1+q^{2nm+nt-2\ell}}}{\lrb{1-q^{nm}}^2\lrb{1-q^{nm+\ell-n}}^2\lrb{1-q^{nm-\ell}}^2} \\
				&= \sum_{k=0}^\infty (-1)^k A_{\ell, n,k}(q) z^{-k}\lrb{1-z}^{2k} = \sum_{k=0}^\infty A_{\ell, n,k}(q) \sum_{j=0}^{2k} \binom{2k}{j} (-z)^{j-k} \\
				&=  \sum_{j=-\infty}^{\infty}\sum_{k=|j|}^\infty  \binom{2k}{j+k} A_{\ell, n,k}(q) (-z)^{j}.
	\end{align*}
	Changing $z$ to $-z$ and then comparing the coefficient of $z^j$ gives us the theorem.
\end{proof}
Next, we prove Theorem \ref{JTP}.
\begin{proof}[Proof of Theorem \ref{JTP}]
We recall the Jacobi triple product identity (see \cite[(2.2.1)]{Hirschhorn})
$$
\sum_{n=-\infty}^{\infty}(-1)^n q^{\frac{n(n+1)}{2}}z^{2n+1}=(z-z^{-1})\prod_{n=1}^{\infty}(1-q^{n})(1-z^{-2}q^{n})(1-z^2q^{n}).
$$
By factoring out $(q;q)_{\infty}$ and letting $z\rightarrow q^\beta z$ gives
\begin{displaymath}
	\begin{split}
		\sum_{n=-\infty}^{\infty}(-1)^nq^{\frac{n(n+1)}{2}+\beta(2n+1)}z^{2n+1}=(q;q)_{\infty}\left(q^{\beta}z-q^{-\beta}z^{-1}\right)\prod_{n=1}^{\infty}(1-z^{-2}q^{n-2\beta})(1-z^2q^{n+2\beta}).
	\end{split}
\end{displaymath}
By a simple re-indexing, we obtain
\begin{displaymath}
	\begin{split}
		\sum_{n=-\infty}^{\infty}(-1)^nq^{\frac{n(n+1)}{2}+\beta(2n+1)}z^{2n+1}=q^{\beta}z(q;q)_{\infty}\prod_{n=0}^{4\beta}\left(1-z^{-2}q^{n-2\beta}\right)\prod_{n=1}^{\infty}(1-z^{-2}q^{n+2\beta})(1-z^2q^{n+2\beta}).
	\end{split}
\end{displaymath}
After straightforward algebraic manipulation, we find
\begin{displaymath}
	\begin{split}
		\sum_{n=-\infty}^{\infty}(-1)^nq^{\frac{n(n+1)}{2}+2n\beta}z^{2n}=(q;q)_{\infty}\prod_{n=0}^{4\beta}\left(1-z^{-2}q^{n-2\beta}\right)\prod_{n=1}^{\infty}\left((1-q^{n+2\beta})^2+(i(z-z^{-1}))^2q^{n+2\beta}\right).
	\end{split}
\end{displaymath}
After factoring out $(q^{2\beta+1};q)_{\infty}^2$ from the infinite product and re-indexing the finite product, we obtain
\begin{displaymath}
	\begin{split}
		\sum_{n=-\infty}^{\infty}(-1)^nq^{\frac{n(n+1)}{2}+2n\beta}z^{2n}=(q;q)_{\infty}(q^{2\beta+1};q)_{\infty}^2\prod_{n=-2\beta}^{2\beta}\left(1-z^{-2}q^{n}\right)\prod_{n=1}^{\infty}\left(1+\frac{q^{n+2\beta}(i(z-z^{-1}))^2}{(1-q^{n+2\beta})^2}\right).
	\end{split}
\end{displaymath}
By definition \eqref{B_k}, we find that the infinite product on the right, as a power series in $(i(z+z^{-1}))^2,$ is the generating function for  $B_{k,\beta}(q)$.   
Therefore, we find that
\begin{align*}
	\sum_{n=-\infty}^{\infty}\frac{(-1)^nq^{\frac{n(n+1)}{2}+2n\beta}}{(q;q)_{\infty}(q^{2\beta+1};q)_{\infty}^2}\cdot z^{2n}=\prod_{n=-2\beta}^{2\beta}\left(1-z^{-2}q^{n}\right)\sum_{n=0}^{\infty}B_{n,\beta}(q)(i(z+z^{-1}))^{2n}.
\end{align*}
Using the Binomial Theorem, followed by a shift in the index of summation, and culminating with a change in the order of summation, we obtain
\begin{align*}
	\sum_{n=-\infty}^{\infty}\frac{(-1)^nq^{\frac{n(n+1)}{2}+2n\beta}}{(q;q)_{\infty}(q^{2\beta+1};q)_{\infty}^2}\cdot z^{2n}&=\prod_{n=-2\beta}^{2\beta}\left(1-z^{-2}q^{n}\right)\sum_{n=0}^{\infty}B_{n,\beta}(q)\sum_{j=0}^{2n}(-1)^{n+j}\binom{2n}{j}z^{2j-2n}\\
	&=\prod_{n=-2\beta}^{2\beta}\left(1-z^{-2}q^{n}\right)\sum_{j=-\infty}^{\infty}\sum_{n=|j|}^{\infty}(-1)^j\binom{2n}{j+n}B_{n,\beta}(q)z^{2j}.
\end{align*}
By replacing $z^2\rightarrow -z$ and using definition \eqref{GF}, we get
\begin{align*}
	\sum_{n=-\infty}^{\infty}\frac{q^{\frac{n(n+1)}{2}+2n\beta+\beta(2\beta+1)}}{(q;q)_{\infty}(q^{2\beta+1};q)_{\infty}^2}\cdot z^{n}=&(1+z^{-1})\left(\sum_{m=-2\beta}^{2\beta}\sum_{t=0}^{2\beta(2\beta+1)}P\left([2\beta],t,m\right)q^tz^m\right)\\
	&\times\left(\sum_{j=-\infty}^{\infty}\sum_{n=|j+2\beta|}^{\infty}\binom{2n}{j+n+2\beta}B_{n,\beta}(q)z^{j}\right).
\end{align*}
After multiplying through $(1+z^{-1})$, we obtain
\begin{displaymath}
\begin{split}
\sum_{n=-\infty}^{\infty}\frac{q^{\frac{n(n+1)}{2}+2n\beta+\beta(2\beta+1)}}{(q;q)_{\infty}(q^{2\beta+1};q)_{\infty}^2}\cdot z^{n}=&\sum_{k=-\infty}^{\infty}\left(\sum_{m+j=k\atop m\in\{-2\beta,\ldots, 2\beta\},j\in\Z}\sum_{t=0}^{2\beta(2\beta+1)} P([2\beta],t,m)q^t\sum_{n=|j+2\beta|}^\infty \binom{2n}{n+j+2\beta}B_{n,\beta}(q)\right.\\
&\left.+\sum_{m+j=k-1\atop m\in\{-2\beta,\ldots, 2\beta\},j\in\Z}\sum_{t=0}^{2\beta(2\beta+1)} P([2\beta],t,m)q^t\sum_{n=|j+2\beta|}^\infty \binom{2n}{n+j+2\beta}B_{n,\beta}(q)\right)z^{k}.
\end{split}
\end{displaymath}
The first part of the theorem follows by comparing the coefficient of $z^{k}$ on both sides. For the proof of second part, we recall the Jacobi triple product identity (see \cite[Theorem~2.8]{Andrews})
$$
\sum_{n=-\infty}^{\infty} q^{n^2}z^{n}=\prod_{n=0}^{\infty}(1-q^{2n+2})(1+z^{-1}q^{2n+1})(1+zq^{2n+1}).
$$
The rest of the proof proceeds similarly to the first part. This completes the proof of the theorem.
\end{proof}


\begin{thebibliography}{99}
	
   \bibitem{AOS} T. ~Amdeberhan, K. Ono, A. Singh,
   \emph{ MacMahon's sums-of-divisors and allied $q$-series},
    Advances in Mathematics (to appear) {\tt https://doi.org/10.48550/arXiv.2311.07496}, (2024).
	
	\bibitem{Andrews} G. E. Andrews, 
	\emph{The theory of partitions}, Cambridge Univ. Press, Cambridge, (1998).
	
	\bibitem{Andrews-Berndt} G. E. Andrews, B. C. Berndt,
	\emph{ Ramanujan's lost notebook, Part III}, Springer, New York, (2012).
	
	\bibitem{Andrews-Rose} G. E. ~Andrews, S. C. F. ~Rose,
	\emph{MacMahon's sum-of-divisors functions, Chebyshev polynomials, and quasimodular forms}, J. Reine Angew. Math. \textbf{676}, 97--103, (2013).
	
	\bibitem{Bachmann} H. Bachmann, \emph{MacMahon's sums-of-divisors and their connections to multiple Eisenstein series}, Res. Number Theory 10, 50, (2024).
	
	\bibitem{R1} H. ~Bachmann and U. ~K\"uhn,
	\emph{The algebra of generating functions for multiple divisor sums and applications to multiple zeta values}, The Ramanujan J., \textbf{40} 3, 605-648, (2016).
	
	\bibitem{Bringmann} K. Bringmann, W. Craig, J.-W. van Ittersum, B. V. Pandey,
	\emph{Limiting behaviour and modular completions of MacMahon-like $q$-series}, {\tt https://arxiv.org/abs/2402.08340}, (2023).
	
	\bibitem{CVO2024} W. Craig, J-W. Van-Ittersum, K. Ono, \emph{Interger partitions detect the primes}, preprint available {\tt https://arxiv.org/abs/2405.06451} (2024).
	
	\bibitem{Gordon} B. Gordon,
	\emph{Some identities in combinatorial analysis}, Quart. J. Math. Oxford
	Ser. \textbf{(2)} 12, 285--290, (1961).
	
	\bibitem{R2} Kh.~Hessami Pilehrood, T.~Hessami Pilehrood,
	\emph{On $q$-analogues of two-one formulas for multiple harmonic sums and multiple zeta star values}, 
	Monatsh. Math. \textbf{176},  275--291, (2015).
	
	\bibitem{Hirschhorn} M. D. Hirschhorn
	\emph{The power of $q$-a personal journey}, Developments in Mathematics, Springer International Publishing AG 2017, 1389-2177, (2018).
	
	\bibitem{Hirschhorn1} M. D. Hirschhorn,
	\emph{A simple proof of an identity of Ramanujan}, J. Aust. Math. Soc. Ser. A \textbf{34}, 31--35, (1983).
	
	\bibitem{R3} D. ~Kreimer, 
	\emph{Knots and Feynman diagrams}, 
	Cambridge Lecture Notes in Physics, vol. 13, Cambridge University Press, Cambridge, (2000). 
	
	\bibitem{R4} T. Q. Th. ~Le and J. ~Murakami, 
	\emph{Kontsevich’s integral for the Homfly polynomial and relations between values of multiple zeta functions}, 
	Topology Appl. 62, no. 2, 193–206, (1995).
	
	\bibitem{MacMahon} P. A. ~MacMahon, 
	\emph{Divisors of Numbers and their Continuations in the Theory of Partitions},
	Proc. London Math. Soc. (2) \textbf{19}, no.1, 75-113 , (1920)
	[also in Percy Alexander MacMahon Collected Papers, Vol.2, pp. 303--341, (ed. G.E. Andrews), MIT Press, Cambridge, (1986)].
	
	\bibitem{OS} K. Ono, A. Singh, \emph{ Remarks on MacMahon’s q-series},
	Journal of Combinatorial Theory, Series A, Vol.207, (2024).
	
	\bibitem{rose} S. C. F. ~Rose,
	\emph{Quasimodularity of generalized sum-of-divisors functions},
	Research in Number Theory \textbf{1}, Paper No. 18, 11 pp, (2015).
	
	\bibitem{Watson} G. N. Watson, 
	\emph{Theorems stated by Ramanujan. $VII$: Theorems on continued
	fractions}, J. London Math. Soc., vol. 4, p. 39--48, (1929).
	
	
\end{thebibliography}
\end{document}